\documentclass[11pt,twoside,leqno]{article}
\usepackage{amsfonts,amssymb,amscd,amsmath,amsthm,enumerate,verbatim,calc,latexsym}
\usepackage[linkcolor=blue,citecolor=magenta,anchorcolor=blue,pdftex=true,bookmarks=true,colorlinks=true]{hyperref}
\newtheorem{thm}{Theorem}[section]
\newtheorem{cor}[thm]{Corollary}
\newtheorem{lem}[thm]{Lemma}

\newtheorem{defn}[thm]{Definition}

\newtheorem{exam}[thm]{Example}
\newtheorem{rem}[thm]{Remark}

\DeclareMathOperator{\Ann}{Ann} 
\DeclareMathOperator{\Ext}{Ext} \DeclareMathOperator{\Supp}{Supp}
\DeclareMathOperator{\Ideals}{Ideals} \DeclareMathOperator{\V}{V}
\DeclareMathOperator{\Hom}{Hom} 
 
\DeclareMathOperator{\depth}{depth} \DeclareMathOperator{\cd}{cd}
 \DeclareMathOperator{\Ass}{Ass}
 
 \DeclareMathOperator{\lc}{H}
 
\DeclareMathOperator{\Spec}{Spec} \DeclareMathOperator{\G}{\Gamma}

\newcommand{\fa}{\mathfrak{a}}
\newcommand{\fb}{\mathfrak{b}}

\newcommand{\fm}{\mathfrak{m}}
\newcommand{\fp}{\mathfrak{p}}

\newcommand{\lo}{\longrightarrow}

\numberwithin{equation}{section}
\begin{document}
\title{\bf Some results on generalized local cohomology modules}
\author{
\\{\bf Alireza Vahidi}\\
\small Department of Mathematics, Payame Noor University, I. R of IRAN\\
\small E-mail: vahidi.ar@pnu.ac.ir
\and
{\bf Moharram Aghapournahr}\\
\small Department of Mathematics, Faculty of Science,  Arak University, Arak, 38156-8-8349, IRAN\\
\small E-mail: m-aghapour@araku.ac.ir
}

\date{}
\maketitle

\begin{abstract}

Let $R$ be a commutative Noetherian ring with non-zero identity, $\fa$ an ideal of $R$, $M$ a finite $R$--module and $X$ an arbitrary $R$--module. Here, we show that, in the Serre subcategories of the category of $R$--modules, how the generalized local cohomology modules, the ordinary local cohomology modules and the extension modules behave similarly at the initial points. We conclude some Artinianness and cofiniteness results for $\lc^{n}_{\fa}(M, X)$, and some finiteness results for $\Supp_R(\lc^{n}_{\fa}(M, X))$ and $\Ass_R(\lc^{n}_{\fa}(M, X))$.\\



\end{abstract}

\section{Introduction}

Let $R$ be a commutative Noetherian ring with non-zero identity. We use symbols $\fa$, $M$, and $X$ as an ideal of $R$, a finite (i.e. finitely generated) $R$--module, and an arbitrary $R$--module which is not necessarily finite. For basic results, notations and terminologies not given in this paper, the reader is referred to \cite {BSh} and \cite {BH}.

The $i$th generalized local cohomology module\\
\centerline{$\lc^i_\fa(M, X)\cong \underset{n\in \mathbb{N}}\varinjlim \Ext^{i}_{R}(M/{\fa}^{n}M, X),$}\\
which is a generalization of the $i$th ordinary local cohomology module\\
\centerline{$\lc^i_\fa(X)\cong \underset{n\in \mathbb{N}}\varinjlim \Ext^{i}_{R}(R/{\fa}^{n}, X),$}\\
was introduced by Herzog in his habilitation \cite {He} and then continued by Suzuki \cite {S}, Bijan-Zadeh \cite {BZ}, Yassemi \cite {Ya} and some other authors. They studied some basic duality theorems, vanishing and other properties of generalized local cohomology modules which also generalize several known facts about extension modules and ordinary local cohomology modules.

In Section \ref {2}, we present the main results of this paper which determine that, for non-negative integers $m$ and $n$, when $R$--modules $\lc^{n}_{\fa}(M, X)$, $\Hom_R(R/\fa, \lc^{n}_{\fa}(M, X))$ and $\Ext^{m}_{R}(M, \lc^{n}_{\fa}(X))$ are in a Serre subcategory of the category of $R$--modules (i.e. the class of $R$--modules which is closed under taking submodules, quotients and extensions), and when $\lc^{m+ n}_{\fa}(M, X)\cong \Ext^{m}_{R}(M, \lc^{n}_{\fa}(X))$ holds (Theorems \ref {2---6}, \ref {2---10}, \ref {2---13} and \ref {2---21}). We use these theorems to show that, in the Serre subcategories of the category of $R$--modules, how the generalized local cohomology modules and the ordinary local cohomology modules behave similarly at the initial points (Corollaries \ref {2---9}, \ref {2---14} and \ref {2---15}). We also find, in Corollaries \ref {2---17} and \ref {2---18}, the relation of regular sequences with respect to Serre classes, introduced  in \cite [Definition 2.6]{AM}, and the membership of generalized local cohomology modules, extension modules and Koszul cohomology modules in Serre subcategories. Note that, one can apply our results to the Serre subcategories of Examples \ref{2---2} and \ref{2---8} to deduce more properties of generalized local cohomology modules.

Section \ref {3} consists of applications. We first, in Corollaries \ref {3---1} and \ref {3---2}, study Artinian generalized local cohomology modules and show that if $\dim_R(R/\fa)= 0$, then the generalized local cohomology modules $\lc^n_\fa(M, X)$ are Artinian and $(\fa+ \Ann_R M)$--cofinite (i.e. $\Supp_R(\lc^n_\fa(M, X))\subseteq \V(\fa+ \Ann_R M)$ and $\Ext_R^i(R/\fa+ \Ann_R M, \lc^n_\fa(M, X))$ is finite for all $i$). Then we present the relation between length, annihilator and support of generalized local cohomology modules, and those of ordinary local cohomology modules (Corollaries \ref {3---3}, \ref {3---4} and \ref {3---5}). We also prove that\\
\centerline{$\underset{i< n}\bigcup \Supp_R(\lc^{i}_{\fa}(M, X))= \underset{i< n}\bigcup \Supp_R(\lc^{i}_{\fa+ \Ann_R M}(X))= \underset{i< n}\bigcup \Supp_R(\Ext^{i}_{R}(M/{\fa}M, X)),$}\\
and if $\lc^{i}_{\fa}(M, X)= 0$ for all $i< n$, then $\Ass_R(\lc^{n}_{\fa}(M, X))= \Ass_R(\Ext^n_{R}(M/{\fa}M, X))$ (Corollaries \ref {3---6} and \ref {3---8}). This implies that if $\Supp_R(\lc^{i}_{\fa}(M, X))$ is finite for all $i< n$, then the finiteness of $\Ass_R(\lc^{n}_{\fa}(M, X))$ is equivalent to the finiteness of $\Ass_R(\Ext^{n}_{R}(M/{\fa}M, X))$. Finally, in the study of finiteness of the set of associated prime ideals of generalized local cohomology modules, we point out the proof of \cite [Theorem 2.3] {Mafi1} contains a flaw, but we show that the statements of \cite [Corollaries 2.4 trough 2.7] {Mafi1} are true (Remark \ref {3---7} and Corollaries \ref {3---11} trough \ref {3---14}).

Even though we can show some of our results by using spectral sequences, we are avoiding the use of this technique completely in this work and we provide more elementary proofs for the results.


\section{Main Results} \label {2}




Let $M$ be a finite $R$--module, $X$ be an arbitrary $R$--module and $n$ be a non-negative integer. We first present sufficient conditions which convince us the $R$--modules $\lc^{n}_{\fa}(M ,X)$ and $\Hom_R(R/\fa, \lc^{n}_{\fa}(M ,X))$ are in a Serre subcategory of the category of $R$--modules.

\begin{defn} \label {2---1}
\emph{Recall that a \textit{Serre subcategory} $\mathcal{S}$ of the category of $R$--modules is a subclass of $R$--modules such that for any short exact sequence\\
\centerline{$0\longrightarrow X'\longrightarrow X\longrightarrow X''\longrightarrow 0,$}\\
the module $X$ is in $\mathcal{S}$ if and only if $X'$ and $X''$ are in $\mathcal{S}$.}
\end{defn}



\begin{exam} \label {2---2}
\emph{The following classes are Serre subcategories of the category of $R$--modules.
\vspace{-2mm}
\begin{itemize}\setlength{\itemsep}{-2mm}
       \item[(a)]  The class of zero $R$--modules.
       \item[(b)]  The class of finite length $R$--modules.
       \item[(c)]  The class of finite $R$--modules.
       \item[(d)]  The class of Artinian $R$--modules.
       \item[(e)]  The class of $R$--modules with finite support.
       \item[(f)]  The class of $R$--modules with Krull dimension less than $n$, where $n$ is a non-negative integer.
       \item[(g)]  The class of $R$--modules with finite Krull dimension.
       \item[(h)]  The class of minimax $R$--modules (An $R$--module $X$ is said to be {\it minimax} if there is a finite submodule $X'$ of $X$ such that $X/X'$ is Artinian \cite {Z}).
\end{itemize}}
\end{exam}



\begin{defn} \label {2---3}
\emph{Let $\lambda: \mathcal{S}\longrightarrow \mathcal{T}$ be a function from a Serre subcategory of the category of $R$--modules $\mathcal{S}$ to a partially ordered Abelian monoid  $(\mathcal{T}, \displaystyle\bigstar, \preceq)$. We say that $\lambda$ is a \textit{subadditive function} if $\lambda(0)= 0$ and for any short exact sequence\\
\centerline{$0\longrightarrow X'\longrightarrow X\longrightarrow X''\longrightarrow 0,$}\\
in which all the terms belong to $\mathcal{S}$, $\lambda(X')\preceq \lambda(X), \lambda(X'')\preceq \lambda(X)$ and $\lambda(X)\preceq \lambda(X')\displaystyle\bigstar \lambda(X'').$}
\end{defn}



\begin{exam} \label {2---4}
\emph{The following functions are subadditive.
\vspace{-2mm}
\begin{itemize}\setlength{\itemsep}{-2mm}
       \item[(a)]  The function $\lambda(X)= l_R(X)$, length of $X$, from the class of finite length $R$--modules to the partially ordered Abelian monoid $(\mathbb{Z},+ ,\leq)$.
       \item[(b)]  The function $\lambda(X)= (0:_R X)$, annihilator of $X$, from the category of $R$--modules to the partially ordered Abelian monoid $(\Ideals(R),. ,\supseteq)$.
       \item[(c)]  The function $\lambda(X)= \Supp_R(X)$, support of $X$, from the category of $R$--modules to the partially ordered Abelian monoid $(\mathrm{P}(\Spec R),\cup ,\subseteq)$.
\end{itemize}}
\end{exam}



In this paper, $\mathcal{S}$ is a Serre subcategory of the category of $R$--modules, $(\mathcal{T}, \displaystyle\bigstar, \preceq)$ is a partially ordered Abelian monoid and $\lambda: \mathcal{S}\longrightarrow \mathcal{T}$ is a subadditive function.

Our method to prove the main results of the paper is based on the induction argument and we need the following useful lemmas for the base cases and inductive steps. Note that, for all $i$, we have the isomorphism\\
\centerline{$\lc^{i}_{\fa}(M, X)\cong \lc^{i}(\G_{\fa}(\Hom_{R}(M, E^\bullet))),$}\\
where $E^\bullet$ is an injective resolution of $X$.

\begin{lem} \label {2---5}
Let $M$ be a finite $R$--module, $X$ be an arbitrary $R$--module and $\fp$ be a prime ideal of $R$. Then the following statements hold true.
\vspace{-2mm}
\begin{itemize}\setlength{\itemsep}{-2mm}
       \item[\emph{(a)}] $\Gamma_{\fa}(M, X)\cong \Hom_{R}(M, \Gamma_\fa (X)).$
       \item[\emph{(b)}] $\lc_{\fa}^i(M, X)_\fp\cong \lc_{{\fa}R_\fp}^i(M_\fp, X_\fp)$ for all $i$.
       \item[\emph{(c)}] If $\Supp_R(M)\cap \Supp_R(X)\subseteq \V(\fa)$, then $\lc^{i}_{\fa}(M, X)\cong\Ext^i_{R}(M, X)$ for all $i$.
\end{itemize}
\end{lem}

\begin{proof}
This is easy and left to the reader.
\end{proof}



\begin{thm} \label {2---6}
Let $M$ be a finite $R$--module, $X$ be an arbitrary $R$--module, and $n$ be a non-negative integer such that $\Ext^{n- r}_{R}(M, \lc^{r}_{\fa}(X))$ is in $\mathcal{S}$ for all $r$, $0\leqslant r\leqslant n.$ Then $\lc^{n}_{\fa}(M, X)\in \mathcal{S}$, and\\
\centerline{$\lambda(\lc^{n}_{\fa}(M, X))\preceq \underset{r= 0}{\overset{n}\bigstar} \lambda(\Ext^{n- r}_{R}(M, \lc^{r}_{\fa}(X))).$}
\end{thm}

\begin{proof}
We prove by using induction on $n$. The case $n=0$ is clear from Lemma \ref {2---5} (a). Suppose that $n> 0$ and that $n- 1$ is settled. Let $\overline{X}= X/\Gamma_{\frak{a}}(X)$ and $L= E(\overline{X})/\overline{X}$ where $E(\overline{X})$ is an injective hull of $\overline{X}$. Since $\Gamma_{\fa}(\overline X)= 0= \Gamma_{\fa}(E(\overline{X}))$, $\Gamma_{\fa}(M, \overline X)= 0= \Gamma_{\fa}(M, E(\overline{X}))$ by Lemma \ref {2---5} (a). Applying the derived functors of $\Gamma_{\frak{a}}(-)$ and $\Gamma_{\frak{a}}(M,-)$ to the short exact sequence\\
\centerline{$0\rightarrow \overline{X}\rightarrow E(\overline{X})\rightarrow L\rightarrow 0,$}\\
we obtain, for all $i> 0$, the isomorphisms\\
\centerline{$\lc^{i-1}_{\fa}(L)\cong \lc^{i}_{\fa}(\overline X)\ (\cong \lc^{i}_{\fa}(X)) \ \ \textmd{and} \ \ \lc^{i-1}_\fa(M, L)\cong \lc^i_\fa(M,\overline X).$}\\
From the above isomorphisms, for all $r$, $0\leqslant r\leqslant n- 1,$ we have\\
\centerline{$\Ext^{(n- 1)- r}_{R}(M,\lc^{r}_{\fa}(L))\cong \Ext^{n- (r+ 1)}_{R}(M, \lc^{r+ 1}_{\fa}( X))$}\\
which is in $\mathcal{S}$ by assumptions. Thus, from the induction hypothesis on $L$,\\
\centerline{$\lc^{n-1}_\fa(M, L)\in \mathcal{S}\ \ $ and $\ \ \lambda(\lc^{n- 1}_{\fa}(M, L))\preceq \underset{r= 0}{\overset{n- 1}\bigstar}\lambda(\Ext^{(n- 1)- r}_{R}(M, \lc^{r}_{\fa}(L)))$.}\\
Therefore\\
\centerline{$\lc^{n}_\fa(M, \overline X)\in \mathcal{S}\ \ $ and $\ \ \lambda(\lc^{n}_{\fa}(M, \overline X))\preceq \underset{r= 1}{\overset{n}\bigstar}\lambda(\Ext^{n- r}_{R}(M, \lc^{r}_{\fa}(X)))$.}\\
Now, by the short exact sequence\\
\centerline{$0\rightarrow \Gamma_{\fa}(X)\rightarrow X\rightarrow \overline X\rightarrow 0$}\\
and Lemma \ref {2---5} (c), we get the long exact sequence\\
\centerline{$\cdots\lo \Ext^{n}_{R}(M, \Gamma_{\fa}(X))\lo \lc^{n}_{\fa}(M, X)\lo \lc^{n}_{\fa}(M, \overline X)\lo \cdots$}\\
which shows that\\
\centerline{$\lc^{n}_\fa(M, X)\in \mathcal{S}\ \ $ and $\ \ \lambda(\lc^{n}_{\fa}(M, X))\preceq \underset{r= 0}{\overset{n}\bigstar}\lambda(\Ext^{n- r}_{R}(M, \lc^{r}_{\fa}(X)))$}\\
as we desired.
\end{proof}



\begin{defn} \label {2---7}
\emph{(\cite [Definition 2.1]{AM} and \cite [Definition 3.1]{ATV}) Recall that, a Serre subcategory of the category of $R$--modules $\mathcal{M}$ is said to be {\it Melkersson subcategory with respect to the ideal $\fa$} if for any $\fa$--torsion $R$--module $X$, $0:_{X}\fa$ is in $\mathcal{M}$ implies that $X$ is in $\mathcal{M}$. Also, $\mathcal{M}$ is called {\it Melkersson subcategory} when it is Melkersson with respect to all ideals of $R$.}
\end{defn}



\begin{exam} \label {2---8}
\emph{The following classes of modules are Melkersson subcategories by Example \ref {2---2} and \cite [Lemma 2.2]{AM}.
\vspace{-2mm}
\begin{itemize}\setlength{\itemsep}{-2mm}
       \item[(a)]  The class of zero $R$--modules.
       \item[(b)]  The class of Artinian $R$--modules.
       \item[(c)]  The class of $R$--modules with finite support.
       \item[(d)]  The class of $R$--modules with Krull dimension less than $n$, where $n$ is a non-negative integer.
       \item[(e)]  The class of $R$--modules with finite Krull dimension.
\end{itemize}}
\end{exam}



In this paper, $\mathcal{M}_\fa$ stands as a Melkersson subcategory with respect to the ideal $\fa$, $\mathcal{M}_{\fa+ \Ann_R M}$ as a Melkersson subcategory with respect to the ideal $\fa+ \Ann_R M$ and $\mathcal{M}$ as a Melkersson subcategory.

The second author and Melkersson in \cite [Theorem 2.9 (i) $\leftrightarrow$ (vi)] {AM} proved the following corollary for Melkersson subcategories, while it was a simple conclusion of Theorem \ref {2---6} for any arbitrary Serre subcategories. This also generalizes \cite [Theorem 2.2 and Corollary 2.3] {sar} for an arbitrary $R$--module $X$ when we consider $\mathcal{S}$ as the class of minimax $R$--modules and the class of Artinian $R$--modules, respectively.

\begin{cor} \label {2---9}
Suppose that $X$ is an arbitrary $R$--module and that $n$ is a non-negative integer. Then the following statements are equivalent.
\vspace{-2mm}
  \begin{itemize}\setlength{\itemsep}{-2mm}
       \item[\emph{(i)}] {$\lc^{i}_{\fa}(X)$ is in $\mathcal S$ for all $i\leqslant n$ \rm{(}for all $i$\rm{)}.}
       \item[\emph{(ii)}] {$\lc^{i}_{\fa}(M, X)$ is in $\mathcal S$ for any finite $R$--module $M$ and for all $i\leqslant n$ \rm{(}for all $i$\rm{)}.}
  \end{itemize}
\end{cor}

\begin{proof}
(i) $\Rightarrow$ (ii). Assume that $i$ is an integer such that $i\leqslant n$. Since $\lc^{r}_{\fa}(X)$ is in $\mathcal S$ for all $r$, $0\leqslant r\leqslant i$, $\Ext^{i- r}_{R}(M, \lc^{r}_{\fa}(X))$ is in $\mathcal{S}$ for all $r$, $0\leqslant r\leqslant i.$ Thus, by Theorem \ref {2---6}, $\lc^{i}_{\fa}(M, X)$ is in $\mathcal{S}.$
\end{proof}



\begin{thm} \label {2---10}
Let $M$ be a finite $R$--module, $X$ be an arbitrary $R$--module, and $n$ be a non-negative integer. Then the following statements hold true.
\vspace{-2mm}
\begin{itemize}\setlength{\itemsep}{-2mm}
       \item[\emph{(a)}] If $\lc^{r}_{\fa}(X)\in \mathcal{S}$ for all $r$, $0\leqslant r< n,$ then $\Hom_R(R/\fa, \lc^{n}_{\fa}(M ,X))\in \mathcal{S}$ whenever $\Ext_R^n(M/\fa M, X)\in \mathcal{S}$.
       \item[\emph{(b)}] If $\lc^{r}_{\fa}(X)= 0$ for all $r$, $0\leqslant r< n,$ then $\Hom_R(R/\fa, \lc^{n}_{\fa}(M ,X))\cong \Ext_R^n(M/\fa M, X)$.
 \end{itemize}
\end{thm}

\begin{proof}
We prove by using induction on $n$.  From Lemma \ref {2---5} (a), we get\\
\centerline{$\begin{array}{llll}
\Hom_R(R/\fa, \Gamma_\fa(M ,X))\!\!
&\cong \ \ \Hom_R(R/\fa, \Hom_R(M, \Gamma_\fa(X)))\\
&\cong \ \ \Hom_R(R/\fa\otimes_R M, \Gamma_\fa(X))\\
&\cong \ \ \Hom_R(M/\fa M, \Gamma_\fa(X))\\
&\cong \ \ \Hom_R(M/\fa M, X)
\end{array}$}\\
because $\Hom_R(M/\fa M, X/\Gamma_\fa(X))= 0$. Thus the assertion follows in the case that $n=0$. Suppose that $n> 0$ and that $n- 1$ is settled. To complete the induction argument, one can use the short exact sequence\\
\centerline{$0\rightarrow \overline{X}\rightarrow E(\overline{X})\rightarrow L\rightarrow 0$}\\
and employ the induction hypothesis with a similar method as in the proof of Theorem \ref {2---6}.
\end{proof}



\begin{rem} \label {2---11}
\emph{Theorem \ref {2---6}, Corollary \ref {2---9} and Theorem \ref {2---10} can be applied to each Serre subcategory mentioned in Example \ref {2---2} resulting in each case in a number of facts about generalized local cohomology modules. One can also use the Serre subcategories and Melkersson subcategories of Examples \ref {2---2} and \ref {2---8} in the results that follow to deduce more properties of generalized local cohomology modules.}
\end{rem}



As an application of the above theorem, we can state the following corollary.

\begin{cor} \label {2---12}
Let $M$ be a finite $R$--module, $X$ be an arbitrary $R$--module, and $n$ be a non-negative integer such that $\lc^{i}_{\fa}(X)\in \mathcal{M}_\fa$ for all $i< n.$ Then $\lc^{n}_{\fa}(M ,X)\in \mathcal{M}_\fa$ whenever $\Ext_R^n(M/\fa M, X)\in \mathcal{M}_\fa$.
\end{cor}

\begin{proof}
Since $\lc^{n}_{\fa}(M ,X)$ is an $\fa$--torsion $R$--module, the assertion follows from Theorem \ref{2---10} (a).
\end{proof}



Now, for non-negative integers $m$ and $n$, we present sufficient conditions which ensure us the $R$--module $\Ext^{m}_{R}(M, \lc^{n}_{\fa}(X))$ is in a Serre subcategory of the category of $R$--modules.

\begin{thm} \label {2---13}
Let $M$ be a finite $R$--module, $X$ be an arbitrary $R$--module, and $m, n$ be non-negative integers. Assume also that
\vspace{-2mm}
\begin{itemize}\setlength{\itemsep}{-2mm}
       \item[\emph{(i)}]   $\lc^{m+ n}_{\fa}(M, X)$ is in $\mathcal S$,
       \item[\emph{(ii)}]  $\Ext^{m+ 1+ r}_{R}(M, \lc^{n- r}_{\fa}(X))$ is in $\mathcal S$ for all $r$, $1\leq r\leq n,$ and
       \item[\emph{(iii)}] $\Ext^{m- 1- r}_{R}(M, \lc^{n+ r}_{\fa}(X))$ is in $\mathcal S$ for all $r$, $1\leq r\leq m- 1.$
\end{itemize}
Then $\Ext^{m}_{R}(M, \lc^{n}_{\fa}(X))\in \mathcal S$, and\\
\centerline{$\lambda(\Ext^{m}_{R}(M, \lc^{n}_{\fa}(X)))\preceq$}\\
\centerline{$\lambda(\lc^{m+ n}_{\fa}(M, X))\bigstar (\underset{r= 1}{\overset{n}\bigstar}\lambda(\Ext^{m+ 1+ r}_{R}(M, \lc^{n- r}_{\fa}(X))))\bigstar (\underset{r= 1}{\overset{m- 1}\bigstar}\lambda(\Ext^{m- 1- r}_{R}(M, \lc^{n+ r}_{\fa}(X)))).$}
\end{thm}

\begin{proof}
We prove by induction on $n$. Let $n=0$ and set $\overline X= X/\G_{\fa}(X)$. By hypothesis (iii), $\Ext^{(m- 1)- r}_{R}(M,\lc^{r}_{\fa}(\overline X))$ is in $\mathcal S$ for all $r$, $0\leqslant r\leqslant m-1$. Thus, from Theorem \ref {2---6},\\
\centerline{$\lc^{m-1}_{\fa}(M, \overline{X})\in \mathcal S\ \ $ and $\ \ \lambda(\lc^{m-1}_{\fa}(M, \overline{X}))\preceq \underset{r= 1}{\overset{m- 1}\bigstar}\lambda(\Ext^{(m- 1)- r}_{R}(M,\lc^{r}_{\fa}(X))).$}\\
By considering Lemma \ref {2---5} (c) and applying the derived functors of $\Gamma_{\fa}(M, -)$ to the short exact sequence\\
\centerline{$0\rightarrow \Gamma_{\fa}(X)\rightarrow X\rightarrow \overline X\rightarrow 0,$}\\
we obtain the long exact sequence\\
\centerline{$\cdots\lo \lc^{m-1}_{\fa}(M, \overline X)\lo \Ext^{m}_{R}(M,\G_{\fa}(X))\lo \lc^{m}_{\fa}(M, X)\lo \cdots$}\\
which shows that, by hypothesis (i), $\Ext^{m}_{R}(M, \Gamma_{\fa}(X))\in \mathcal{S}$ and\\
\centerline{$\lambda(\Ext^{m}_{R}(M,\G_{\fa}(X)))\preceq \lambda(\lc^{m}_{\fa}(M, X))\bigstar(\underset{r= 1}{\overset{m- 1}\bigstar}\lambda(\Ext^{m- 1- r}_{R}(M,\lc^{r}_{\fa}(X))))$.}\\
Thus the assertion follows in this case.

Now, assume that $n> 0$ and that $n- 1$ is settled. Let $\overline{X}= X/\Gamma_{\frak{a}}(X)$ and $L= E(\overline{X})/\overline{X}$ where $E(\overline{X})$ is an injective hull of $\overline{X}$. By the short exact sequence\\
\centerline{$0\rightarrow \overline{X}\rightarrow E(\overline{X})\rightarrow L\rightarrow 0,$}\\
the proof is sufficiently similar to that of Theorem \ref {2---6} to be omitted. We leave the proof to the reader.
\end{proof}



The next corollary shows that, in Melkersson subcategories, the generalized local cohomology modules $\lc^{i}_{\fa}(M, X)$ and the ordinary local cohomology modules $\lc^{i}_{\fa+ \Ann_R M}(X)$ behave similarly at the initial points.

\begin{cor} \label {2---14}
Let $M$ be a finite $R$--module, $X$ be an arbitrary $R$--module and $n$ be a non-negative integer. Then the following statements are equivalent.
\vspace{-2mm}
\begin{itemize}\setlength{\itemsep}{-2mm}
          \item[\emph{(i)}]  {$\lc^{i}_{\fa+ \Ann_R M}(X)$ is in $\mathcal M_{\fa+ \Ann_R M}$ for all $i\leqslant n$ \rm{(}for all $i$\rm{)}.}
          \item[\emph{(ii)}]  {$\lc^{i}_{\fa}(M, X)$ is in $\mathcal M_{\fa+ \Ann_R M}$ for all $i\leqslant n$ \rm{(}for all $i$\rm{)}.}
\end{itemize}
\end{cor}

\begin{proof}
(i) $\Rightarrow$ (ii). Since $\lc^{i}_{\fa}(M, X)\cong \lc^{i}_{\fa+ \tiny\Ann_R M}(M, X)$ for all $i$, the assertion holds from Corollary \ref {2---9}.

(ii) $\Rightarrow$ (i). We use induction on $n$. Let $n=0$. By considering the exact sequence\\
\centerline{$0\lo \Hom_{R}(M/{\fa}M, \Gamma_{\fa}(X))\lo \Hom_{R}(M, \Gamma_{\fa}(X)),$}\\
$\Hom_{R}(M/{\fa}M, \Gamma_{\fa}(X))$ is in $\mathcal M_{\fa+ \Ann_R M}$ from Lemma \ref {2---5} (a). Thus, by \cite [Theorem 2.9 (iv) $\rightarrow$ (i)] {AM}, $\Gamma_{\fa+ \Ann_R M}(\Gamma_{\fa}(X))$ is in $\mathcal M_{\fa+ \Ann_R M}$. Therefore $\Gamma_{\fa+ \tiny\Ann_R M}(X)$ is in $\mathcal M_{\fa+ \Ann_R M}$.

Assume that $n> 0$ and that $n- 1$ is settled. By the induction hypothesis, $\lc^{i}_{\fa+ \Ann_R M}(X)$ is in $\mathcal M_{\fa+ \Ann_R M}$ for all $i\leqslant n- 1$. Apply Theorem \ref {2---13} with $ m= 0$ to see that $\Hom_{R}(M, \lc^{n}_{\fa+ \tiny\Ann_R M}(X))$ is in $\mathcal M_{\fa+ \Ann_R M}$. Thus $\Hom_{R}(M/{\fa}M, \lc^{n}_{\fa+ \tiny\Ann_R M}(X))$ is in
$\mathcal M_{\fa+ \Ann_R M}$ by the exact sequence\\
\centerline{$0\lo \Hom_{R}(M/{\fa}M, \lc^{n}_{\fa+ \tiny\Ann_R M}(X))\lo \Hom_{R}(M, \lc^{n}_{\fa+ \tiny\Ann_R M}(X)).$}\\
Again from \cite [Theorem 2.9 (iv) $\rightarrow$ (i)] {AM}, we have $\Gamma_{\fa+ \Ann_R M}(\lc^{n}_{\fa+ \tiny\Ann_R M}(X))\in \mathcal M_{\fa+ \Ann_R M}$ which shows that $\lc^{n}_{\fa+ \tiny\Ann_R M}(X)\in \mathcal M_{\fa+ \Ann_R M}$.
\end{proof}



\begin{cor} \label {2---15}
Let $M$ be a finite $R$--module and $X$ be an arbitrary $R$--module. Then we have
\vspace{-2mm}
\begin{itemize}\setlength{\itemsep}{-2mm}
       \item[\emph{(a)}] $\inf \{i: \lc^{i}_{\fa}(M, X)\notin \mathcal{S}\}\geq \inf \{i: \lc^{i}_{\fa}( X)\notin \mathcal{S}\}.$
       \item[\emph{(b)}] $\inf \{i: \lc^{i}_{\fa}(M, X)\notin \mathcal{M}_{\fa+ \Ann_R M}\}= \inf \{i: \lc^{i}_{\fa+ \Ann_R M}( X)\notin \mathcal{M}_{\fa+ \Ann_R M}\}.$
       \item[\emph{(c)}] $\inf \{i: \lc^{i}_{\fa}(M, X)\notin \mathcal{M}_{\fa+ \Ann_R M}\}= \inf \{i: \lc^{i}_{\fa}(X)\notin \mathcal{M}_{\fa+ \Ann_R M}\}$ whenever $\Gamma_{\Ann_R M}(X)= X$.
       \item[\emph{(d)}] $\inf \{i: \lc^{i}_{\fa}(M, X)\notin \mathcal{M}_\fa\}= \inf \{i: \lc^{i}_{\fa}(X)\notin \mathcal{M}_\fa\}$ whenever $\Ann_R M\subseteq \fa$ \rm{(}e.g. $M$ is faithful\rm{)}.
\end{itemize}
\end{cor}

\begin{proof}
Follows from Corollaries \ref {2---9} and \ref {2---14}.
\end{proof}



In the next corollary, we state the membership of the generalized local cohomology modules with respect to different ideals in Melkersson subcategories of the category of $R$--modules.

\begin{cor} \label {2---16}
Suppose that $M$ is a finite $R$--module and $X$ is an arbitrary $R$--module. Assume also that $n$ is a non-negative integer and $\fb$ is an ideal of $R$ such that $\fa\subseteq \fb$. Then $\lc^{i}_{\fb}(M, X)$ is in $\mathcal{M}$ for all $i\leq n$ \rm{(}for all $i$\rm{)} whenever $\lc^{i}_{\fa}(M, X)$ is in $\mathcal{M}$ for all $i\leq n$ \rm{(}for all $i$\rm{)}.
\end{cor}

\begin{proof}
Follows from Corollary \ref {2---14} and \cite [Proposition 3.4] {ATV}.
\end{proof}



In \cite [Definition 2.6 and Example 2.8]{AM}, the second author and Melkersson  introduced the concept of $\mathcal{S}$--regular sequences on a module that recovered poor sequences, filter-regular sequences, generalized regular sequences and sequences in dimension$> n$, where $n$ is a non-negative integer, on a module. They also found, in \cite [Theorem 2.9 (i) $\leftrightarrow$ (vii)]{AM}, the relation of this notion on a finite module and the membership of local cohomology modules in Melkersson subcategories. In the next corollary, we state a similar characterization for generalized local cohomology modules. Coung and Hoang in \cite [Theorem 3.1]{ CH1} proved part [(i) $\leftrightarrow$ (iv)] of the following corollary for the class of Artinian $R$--modules in the case that $R$ was a local ring.

\begin{cor} \label {2---17}
Suppose that $M$ is a finite $R$--module such that $\fa+ \Ann_R M= (x_1,\dots,x_r)$. Assume also that $X$ is an arbitrary $R$--module and that $n$ is a non-negative integer. Then the following statements are equivalent.
\vspace{-2mm}
\begin{itemize}\setlength{\itemsep}{-2mm}
       \item[\emph{(i)}] {$\lc^{i}_{\fa}(M, X)$ is in $\mathcal M_{\fa+ \Ann_R M}$ for all $i\leqslant n$ \rm{(}for all $i$\rm{)}.}
       \item[\emph{(ii)}] {$\Ext^{i}_{R}(M/{\fa}M,X)$ is in $\mathcal M_{\fa+ \Ann_R M}$ for all $i\leqslant n$ \rm{(}for all $i$\rm{)}.}
       \item[\emph{(iii)}] {$\lc^i(x_1, \dots,x_r; X)$ is in $\mathcal M_{\fa+ \Ann_R M}$ for all $i\leqslant n$ \rm{(}for all $i$\rm{)}.}
       \end{itemize}
       When $X$ is finite, these conditions are also equivalent to:
       \vspace{-2mm}
       \begin{itemize}\setlength{\itemsep}{-2mm}
       \item[\emph{(iv)}] {There is a sequence of length $n+ 1$ in $\fa+ \Ann_R M$ that is $\mathcal M_{\fa+ \Ann_R M}$--regular on $X$.}
      \end{itemize}
\end{cor}

\begin{proof}
This follows from Corollary \ref {2---14} and \cite [Theorem 2.9] {AM}.
\end{proof}



Suppose that $X$ is a finite $R$--module such that $X/\fa X$ is not in $\mathcal{M}_\fa$. The second author and Melkersson, in \cite [Lemma 2.14] {AM}, proved that every sequence in $\fa$ which is $\mathcal{M}_\fa$--regular on $X$ can be extended to a maximal one and all maximal $\mathcal{M}_\fa$--regular sequences on $X$ in $\fa$ have the same length. They denoted this common length by $\mathcal{M}_\fa$--$\depth_{\fa}(X)$, in \cite [Definition 2.15] {AM}, and proved, in \cite [Theorem 2.18] {AM}, that it is the least integer such that $\lc^{i}_{\fa}(X)$, $\Ext^{i}_{R}(R/{\fa},X)$ or Koszul cohomology modules with respect to $\fa$ are not in $\mathcal{M}_\fa$. Using the Melkersson subcategories of Example \ref {2---8}, this notion gives ordinary depth, filter-depth, generalized depth and $n$-depth, where $n$ is a non-negative integer. In the following, we prove that $\mathcal{M}_{\fa+ \Ann_R M}-\depth_{\fa+ \Ann_R M}(X)$ is the least integer such that $\lc^{i}_{\fa}(M, X)$, $\Ext^{i}_{R}(M/{\fa}M, X)$ or Koszul cohomology modules with respect to $\fa+ \Ann_R M$ are not in $\mathcal{M}_{\fa+ \Ann_R M}$. This generalizes the result of Bijan-Zadeh \cite [Proposition 5.5] {BZ} when we consider $\mathcal{M}_{\fa+ \Ann_R M}$ as the class of zero $R$--modules. It also recovers \cite [Theorem 2.2] {CT}, \cite [Theorem 3.1] {CH1}, \cite [Theorem 4.1] {CH} and \cite [Theorem 2.8] {Mafi3} if we put $\mathcal{M}_{\fa+ \Ann_R M}$ the class of Artinian $R$--modules or the class of $R$--modules with finite support. Note that, all of these theorems are in the local case while our corollary is in general.

\begin{cor} \label {2---18}
Suppose that $M$ is a finite $R$--module with $\fa+ \Ann_R M= (x_1,\dots,x_r)$ and  $X$ is a finite $R$--module with $X/(\fa+ \Ann_R M)X\notin \mathcal{M}_{\fa+ \Ann_R M}$. Then
\vspace{-2mm}
\begin{itemize}\setlength{\itemsep}{-2mm}
       \item[\emph{(a)}] $\mathcal{M}_{\fa+ \Ann_R M}-\depth_{\fa+ \Ann_R M}(X)=\inf \{i: \lc^{i}_{\fa}(M, X)\notin \mathcal{M}_{\fa+ \Ann_R M}\}.$
       \item[\emph{(b)}] $\mathcal{M}_{\fa+ \Ann_R M}-\depth_{\fa+ \Ann_R M}(X)=\inf \{i: \Ext^{i}_{R}(M/{\fa}M,X)\notin \mathcal{M}_{\fa+ \Ann_R M}\}.$
       \item[\emph{(c)}] $\mathcal{M}_{\fa+ \Ann_R M}-\depth_{\fa+ \Ann_R M}(X)=\inf \{i: \lc^i(x_1,\dots,x_r ; X)\notin \mathcal{M}_{\fa+ \Ann_R M}\}.$
\end{itemize}
\end{cor}

\begin{proof}
Follows from \cite [Lemma 2.14] {AM} and Corollary \ref {2---17}.
\end{proof}



As applications of Theorems \ref {2---6} and \ref {2---13}, we can state the following corollaries.

\begin{cor} \label {2---19}
Let $M$ be a finite $R$--module, $X$ be an arbitrary $R$--module, and $n$ be a non-negative integer such that $\Ext^{j-i}_{R}(M,\lc^{i}_{\fa}(X))$ is in $\mathcal{S}$ for all $i$, $j$ with $0\leq i\leq n-1$ and $j= n, n+1$. Then $\lc^{n}_{\fa}(M, X)$ is in $\mathcal{S}$ if and only if $\Hom_{R}(M, \lc^{n}_{\fa}(X))$ is in $\mathcal{S}$.
\end{cor}



\begin{cor} \label {2---20}
Let $X$ be an $R$--module and $m, n$ be non-negative integers such that $\lc^{i}_{\fa}(X)$ is in $\mathcal S$ for all $i$, $0\leqslant i \leqslant n-1$ or $n+ 1\leqslant i\leqslant m+ n$. Then $\lc^{m+ n}_{\fa}(M, X)$ is in $\mathcal S$ if and only if $\Ext^{m}_{R}(M, \lc^{n}_{\fa}(X))$ is in $\mathcal S$.
\end{cor}



In the following theorem, for non-negative integers $m$ and $n$, we find some sufficient conditions for validity of the isomorphism $\lc^{m+ n}_{\fa}(M, X)\cong \Ext^{m}_{R}(M, \lc^{n}_{\fa}(X))$.

\begin{thm} \label {2---21}
Let $M$ be a finite $R$--module, $X$ be an arbitrary $R$--module, and $m, n$ be non-negative integers. Assume also that
\vspace{-2mm}
\begin{itemize}\setlength{\itemsep}{-2mm}
      \item[\emph{(i)}] $\Ext^{m+ n- r}_{R}(M, \lc^{r}_{\fa}(X))= 0$ for all $r$, $0\leq r\leq n- 1$ or $n+ 1\leq r\leq m+ n$,
      \item[\emph{(ii)}] $\Ext^{m+ 1+ r}_{R}(M, \lc^{n- r}_{\fa}(X))= 0$ for all $r$, $1\leq r\leq n,$ and
      \item[\emph{(iii)}] $\Ext^{m- 1- r}_{R}(M,\lc^{n+ r}_{\fa}(X))= 0$ for all $r$, $1\leq r\leq m- 1.$
\end{itemize}
Then we have $\lc^{m+ n}_{\fa}(M, X)\cong \Ext^{m}_{R}(M, \lc^{n}_{\fa}(X))$.
\end{thm}

\begin{proof} We prove by using induction on $n$. Let $n= 0$. We have $\lc^{m- 1}_{\fa}(M, X/\Gamma_{\fa}(X))= 0= \lc^{m}_{\fa}(M, X/\Gamma_{\fa}(X))$ from hypothesis (iii) and (i), and Theorem \ref {2---6} with $\mathcal{S}= 0$. Now, the assertion follows by the exact sequence\\
\centerline{$\lc^{m- 1}_{\fa}(M, X/\Gamma_{\fa}(X))\lo \Ext^{m}_{R}(M, \Gamma_{\fa}(X))\lo \lc^{m}_{\fa}(M, X)\lo \lc^{m}_{\fa}(M, X/\Gamma_{\fa}(X))$}\\
obtained from the short exact sequence\\
\centerline{$0\lo \Gamma_{\fa}(X)\lo X\lo X/\Gamma_{\fa}(X)\lo 0$}\\
and Lemma \ref {2---5} (c).

Assume that $n> 0$ and that $n- 1$ is settled. By considering the short exact sequence\\
\centerline{$0\rightarrow \overline{X}\rightarrow E(\overline{X})\rightarrow L\rightarrow 0,$}\\
the proof is similar to that of Theorem \ref {2---6}.
\end{proof}



Yassemi, in \cite [Example 3.6]{Ya}, has given an example to show that the $R$--modules $\lc^{n}_{\fa}(M, X)$ and $\Hom_{R}(M, \lc^{n}_{\fa}(X))$ are not always equal. We show that, with some conditions, they are isomorph.

\begin{cor} \label {2---22} \emph{(cf. \cite [Proposition 2.3 (ii)]{HV})}
Let $M$ be a finite $R$--module, $X$ be an arbitrary $R$--module and $n$ be a non-negative integer such that $\Ext^{j-i}_{R}(M, \lc^{i}_{\fa}(X))= 0$ for all $i$, $j$ with $0\leqslant i\leqslant n-1$ and $j= n, n+ 1$. Then we have $\lc^{n}_{\fa}(M, X)\cong \Hom_{R}(M, \lc^{n}_{\fa}(X))$.
\end{cor}

\begin{proof}
Apply Theorem \ref {2---21} with $m= 0$.
\end{proof}



In consistence with Corollary \ref {2---20}, one can state the following corollary which shows that if $X$ is a finite module and $\fa$ is an ideal generated by an $X$–-regular sequence of length $n$, then the generalized local cohomology modules are exactly extension modules of ordinary local cohomology modules.

\begin{cor} \label {2---23}
Suppose that $M$ is a finite $R$--module, $X$ is an arbitrary $R$--module, and $n, m$ are non-negative integers such that $n\leq m$. Assume also that $\lc^{i}_{\fa}(X)= 0$ for all $i$, $i\neq n$ \rm{(}resp. $0\leq i\leq n-1$ or $n+1\leq i\leq m$\rm{)}. Then we have $\lc^{i+ n}_{\fa}(M, X)\cong \Ext^{i}_{R}(M, \lc^{n}_{\frak{a}}(X))$ for all $i$, $i\geq 0$ \rm{(}resp. $0\leq i\leq m-n$\rm{)}.
\end{cor}

\begin{proof}
For all $i$, $i\geq 0$ {\rm(}resp. $0\leq i\leq m-n${\rm)}, apply Theorem \ref {2---21} with $m= i$.
\end{proof}


\section{Applications} \label {3}




Recall that, an $R$--module $X$ is said to be {\it $\fa$--cofinite} if $\Supp_{R}(X)\subseteq V(\fa)$ and $\Ext^{i}_{R}(R/\fa, X)$ is finite for all $i$. Note that, by \cite [Proposition 4.1] {Mel}, the class of Artinian $\fa$--cofinite modules is a Melkersson subcategory with respect to the ideal $\fa$.

\begin{cor} \label {3---1}
Let $M$ be a finite $R$--module, $X$ be an arbitrary $R$--module and $n$ be a non-negative integer. Then the following statements are equivalent.
\vspace{-2mm}
\begin{itemize}\setlength{\itemsep}{-2mm}
\item[\emph{(i)}]   {$\lc^{i}_{\fa}(M, X)$ is Artinian and $(\fa+ \Ann_R M)$--cofinite for all $i\leqslant n$ \rm{(}for all $i$\rm{)}.}
\item[\emph{(ii)}]  {$\Ext_R^i(M/{\fa}M, X)$ has finite length for all $i\leqslant n$ \rm{(}for all $i$\rm{)}.}
\end{itemize}
\end{cor}

\begin{proof}
(i) $\Rightarrow$ (ii). From Corollary \ref {2---14}, $\lc^{i}_{\fa+ \Ann_R M}(X)$ is Artinian and $(\fa+ \Ann_R M)$--cofinite for all $i\leqslant n$. Thus $\Ext_R^i(R/{\fa+ \Ann_R M}, X)$ has finite length for all $i\leqslant n$ by \cite [Corollary 4.12] {ATV}. Therefore, from \cite [Proposition 3.4] {HV}, $\Ext_R^i(M/{\fa}M, X)$ has finite length for all $i\leqslant n.$

(ii) $\Rightarrow$ (i). Since every finite length $(\fa+ \Ann_R M)$--torsion module is Artinian and $(\fa+ \Ann_R M)$--cofinite, the assertion follows from Corollary \ref {2---17}.
\end{proof}



Chu and Tang in \cite [Proposition 2.4] {CT} proved the part [(i) $\leftrightarrow$ (ii)] of the following corollary in  the local case (see also \cite [Corrollary 3.2] {CH1}, \cite [Corollary 3.3] {DV} and \cite [Theorem 2.2] {DST}).

\begin{cor} \label {3---2}
Suppose that $M, X$ are finite $R$--modules and that $n$ is a non-negative integer. Then the following statements are equivalent.
\vspace{-2mm}
\begin{itemize}\setlength{\itemsep}{-2mm}
\item[\emph{(i)}]   {$\dim_R(\lc^{i}_{\fa}(M, X))\leq 0$ for all $i\leqslant n$ \rm{(}for all $i$\rm{)}.}
\item[\emph{(ii)}]  {$\lc^{i}_{\fa}(M, X)$ is Artinian for all $i\leqslant n$ \rm{(}for all $i$\rm{)}.}
\item[\emph{(iii)}] {$\lc^{i}_{\fa}(M, X)$ is Artinian and $(\fa+ \Ann_R M)$--cofinite for all $i\leqslant n$ \rm{(}for all $i$\rm{)}.}
\end{itemize}
In particular, if $\dim_R(R/\fa)= 0$, then $\lc^{i}_{\fa}(M, X)$ is Artinian and $(\fa+ \Ann_R M)$--cofinite for all $i$.
\end{cor}

\begin{proof}
(i) $\Rightarrow$ (iii). Since every finite module with zero dimension is of finite length, the assertion follows from Corollary \ref {2---17} [(i) $\rightarrow$ (ii)] (where $\mathcal{M}_{\fa+ \Ann_R M}$ is taken the class of $R$--modules with Krull dimension less than $1$) and Corollary \ref {3---1}  [(ii) $\rightarrow$ (i)].
\end{proof}



In \cite [Theorem 3.2] {Sch} and for a non-negative integer $n$, Schenzel proved that
\vspace{-2mm}
\begin{itemize}\setlength{\itemsep}{-2mm}
       \item $\Ext^{n}_{R}(M, X)$ is of finite length, and
       \item  $l_R(\Ext^{n}_{R}(M, X))\leq \displaystyle\sum_{i= 0}^{n}l_R(\Ext^{n- i}_{R}(M, \lc^{i}_{\fm}(X)))$
\end{itemize}
when $(R, \fm)$ is a local ring and $M, X$ are finite $R$--modules such that $M\otimes_R X$ is of finite length. As an application of Theorem \ref {2---6}, by considering Lemma \ref{2---5} (c) and \cite [Lemma 3.1] {Sch}, the following corollary extends \cite [Theorem 3.2] {Sch}.

\begin{cor} \label {3---3}
Let $M$ be a finite $R$--module, $X$ be an arbitrary $R$--module and $n$ be a non-negative integer such that $\Ext^{n- i}_{R}(M, \lc^{i}_{\fa}(X))$ is of finite length for all $i\leqslant n.$ Then
\vspace{-2mm}
\begin{itemize}\setlength{\itemsep}{-2mm}
       \item[\emph{(a)}] $\lc^{n}_{\fa}(M, X)$ is  of finite length, and
       \item[\emph{(b)}] $l_R(\lc^{n}_{\fa}(M, X))\leq \displaystyle\sum_{i= 0}^{n}l_R(\Ext^{n- i}_{R}(M, \lc^{i}_{\fa}(X))).$
\end{itemize}
\end{cor}

\begin{proof}
Since the class of finite length $R$--modules is a Serre subcategory of the category of $R$--modules and $\lambda(X)= l_R(X)$ is a subadditive function from the class of finite length $R$--modules to the partially ordered Abelian monoid $(\mathbb{Z},+ ,\leq)$, the assertion follows form Theorem \ref {2---6}.
\end{proof}



As another application of Theorem \ref {2---6}, we find the relation between annihilator of generalized local cohomology modules and annihilator of ordinary local cohomology modules.

\begin{cor} \label {3---4} Let $M$
be a finite $R$--module, $X$ be an arbitrary $R$--module and $n$ be a non-negative integer. Then we have
\vspace{-2mm}
\begin{itemize}\setlength{\itemsep}{-2mm}
       \item[\emph{(a)}] $\displaystyle\prod_{i= 0}^{n}(0:_R \Ext^{n- i}_{R}(M, \lc^{i}_{\fa}(X)))\subseteq (0:_R \lc^{n}_{\fa}(M, X)).$
       \item[\emph{(b)}] $\displaystyle\prod_{i= 0}^{n}(0:_R\lc^{i}_{\fa}(X))\subseteq \displaystyle\bigcap_{i= 0}^{n}(0:_R \lc^{i}_{\fa}(M, X)).$
\end{itemize}
\end{cor}

\begin{proof}
(a) Since $\lambda(X)= (0:_R X)$ is a subadditive function from the category of $R$--modules to the partially ordered Abelian monoid $(\Ideals(R),. ,\supseteq)$, the assertion follows form Theorem \ref {2---6}.

(b) For all $i\leqslant j\leqslant n$, we have $(0:_R\lc^{i}_{\fa}(X))\subseteq (0:_R \Ext^{j- i}_{R}(M, \lc^{i}_{\fa}(X)))$. Thus the assertion follows from part (a).
\end{proof}



In the course of the remaining parts of the paper for an ideal $\fa$ of $R$ and for an arbitrary $R$--module $X$, by $\cd_R(\fa, X)$ (cohomological dimension of $X$ with respect to $\fa$), we mean the largest integer $i$ in which $\lc^{i}_{\fa}(X)$ is non-zero. The next result presents the relation between support of generalized local cohomology modules and support of ordinary local cohomology modules.

\begin{cor} \label {3---5}
Let $M$ be a finite $R$--module, $X$ be an arbitrary $R$--module and $n$ be a non-negative integer. Then we have
\vspace{-2mm}
\begin{itemize}\setlength{\itemsep}{-2mm}
       \item[\emph{(a)}] $\Supp_R(\lc^{n}_{\fa}(M, X))\subseteq \displaystyle\bigcup_{i\leq n}\Supp_R(\Ext^{n- i}_{R}(M, \lc^{i}_{\fa}(X))).$
       \item[\emph{(b)}] $\displaystyle\bigcup_{i\leq n}\Supp_R(\lc^{i}_{\fa}(M, X))\subseteq \displaystyle\bigcup_{i\leq n}\Supp_R(\lc^{i}_{\fa}(X)).$
\end{itemize}
In particular, $\Supp_R(\lc^{n}_{\fa}(M, X))\subseteq \displaystyle\bigcup_{i\leq \cd_R(\fa, X)}\Supp_R(\lc^{i}_{\fa}(X))$.
\end{cor}

\begin{proof}
(a) Since $\lambda(X)= \Supp_R(X)$ is a subadditive function from the category of $R$--modules to the partially ordered Abelian monoid $(\mathrm{P}(\Spec R),\cup ,\subseteq)$, the assertion follows form Theorem \ref {2---6}.

(b) Follows from the first part.
\end{proof}



The vanishing of generalized local cohomology modules from upper bounds needs special conditions and in all of them $M$ must have finite projective dimension (see \cite [Theorems 2.5 and 3.7] {Ya}, \cite [Theorem  3.1] {CH} and \cite [Proposition 2.8] {HV}). However, in the following corollary, we show that there is a union of finitely many supports of generalized local cohomology modules such that the other supports can be viewed as its subset even if $M$ has infinite projective dimension. Parts (a) and (b) of the following corollary in the local case has been proven in \cite [Lemma 2.8 and Corollary 2.9] {CH} by Coung and Hoang when $X$ is a finite $R$--module but we prove it without assuming that $X$ is finite and with no restrictions on $R$.

\begin{cor} \label {3---6}
Let $M$ be a finite $R$--module, $X$ be an arbitrary $R$--module and $n$ be a non-negative integer. Then the following statements hold
true.
\vspace{-2mm}
\begin{itemize}\setlength{\itemsep}{-2mm}
       \item[\emph{(a)}] $\displaystyle\bigcup_{i\leq n} \Supp_R(\lc^{i}_{\fa}(M,X))= \displaystyle\bigcup_{i\leq n} \Supp_R(\lc^{i}_{\fa+\Ann_R{M}}(X))$.
       \item[\emph{(b)}] $\displaystyle\bigcup_{i\leq n} \Supp_R(\lc^{i}_{\fa}(M,X))= \displaystyle\bigcup_{i\leq n} \Supp_R(\Ext^{i}_{R}(M/{\fa}M,X))$.
       \item[\emph{(c)}] $\displaystyle\bigcup_{i\leq n} \Supp_R(\lc^{i}_{\fa}(M,X))$ is a closed set when $X$ is a finite $R$--module.
       \item[\emph{(d)}] $\displaystyle\bigcup_{i} \Supp_R(\lc^{i}_{\fa}(M,X))$ is a closed set when $X$ is a finite $R$--module.
\end{itemize}
In particular, $\Supp_R(\lc^{n}_{\fa}(M, X))\subseteq \displaystyle\bigcup_{i\leq \cd_R(\fa+ \Ann_R M, X)}\Supp_R(\lc^{i}_{\fa}(M, X))$.
\end{cor}

\begin{proof}
(a) By Lemma \ref {2---5} (b) and Corollary \ref {2---14}, we have\\
\centerline{$\begin{matrix}
\fp\notin \underset{i\leq n}\bigcup \Supp_R(\lc^{i}_{\fa}(M,X))
     &\Leftrightarrow &\forall i\leq n;&\lc^{i}_{\fa}(M,X)_{\fp}= 0\\
     &\Leftrightarrow &\forall i\leq n;&\lc^{i}_{{\fa}R_{\fp}}(M_{\fp},X_{\fp})= 0\\
     &\Leftrightarrow &\forall i\leq n;&\lc^{i}_{{\fa}R_{\fp}+{\Ann_{R_{\fp}}M_{\fp}}}(X_{\fp})= 0\\
     &\Leftrightarrow &\forall i\leq n;&\lc^{i}_{\fa+\Ann_R{M}}(X)_{\fp}= 0\\
     &\Leftrightarrow &&\fp\notin \underset{i\leq n}\bigcup \Supp_R(\lc^{i}_{\fa+\Ann_R{M}}(X))
\end{matrix}$}\\
as we desired.

(b) From Lemma \ref {2---5} (b) and Corollary \ref {2---17} [(i)
$\leftrightarrow$ (ii)], we get\\
\centerline{$\begin{matrix}
\fp\notin \underset{i\leq n}\bigcup \Supp_R(\lc^{i}_{\fa}(M,X))
     &\Leftrightarrow &\forall i\leq n;&\lc^{i}_{\fa}(M,X)_{\fp}= 0\\
     &\Leftrightarrow &\forall i\leq n;&\lc^{i}_{{\fa}R_{\fp}}(M_{\fp},X_{\fp})= 0\\
     &\Leftrightarrow &\forall i\leq n;&\Ext^{i}_{R_{\fp}}(M_{\fp}/({\fa}R_{\fp})M_{\fp},X_{\fp})= 0\\
     &\Leftrightarrow &\forall i\leq n;&\Ext^{i}_{R}(M/{\fa}M, X)_{\fp}= 0\\
     &\Leftrightarrow &&\fp\notin \underset{i\leq n}\bigcup \Supp_R(\Ext^{i}_{R}(M/{\fa}M,X))
\end{matrix}$}\\
as desired.

(c) This is clear from the second part.

(d) By the first part, we have\\
\centerline{$\displaystyle\bigcup_{i} \Supp_R(\lc^{i}_{\fa}(M, X))= \displaystyle\bigcup_{i\leq \cd_R(\fa+ \Ann_R M, X)} \Supp_R(\lc^{i}_{\fa}(M,X)).$}\\
Thus the assertion follows from part (c).
\end{proof}



\begin{rem} \label {3---7}
Let $M$ be a finite $R$--module, $X$ be an arbitrary $R$--module and $n$ be a non-negative integer. In \emph{\cite [Theorem 2.3]{Mafi1}}, Mafi proved\\
\centerline{$\Ass_R(\lc^{n}_{\fa}(M, X))\subseteq \bigcup_{i= 0}^{n}\Ass_R(\Ext^{n- i}_{R}(M, \lc^{i}_{\fa}(X)))$}\\
and, in \emph{\cite [Corollaries 2.4 through 2.7] {Mafi1}}, used it to deduce some results about finiteness of the set of associated prime ideals of generalized local cohomology modules. Although Corollaries \emph{2.4} through \emph{2.7} in \emph{\cite {Mafi1}} are true, the proof of \emph{\cite [Theorem 2.3]{Mafi1}} holds a flaw. In its proof, even though $(E_\infty^{i, t- i}=)\ker d_{t+ 2}^{i, t- i}$ is a subquotient of $\ker d_{2}^{i, t- i}(\subseteq E_2^{i, t- i})$,\\
\centerline{$(E_\infty^{i, t- i}=)\ \ \ker d_{t+ 2}^{i, t- i}\ \subseteq \ \ker d_{2}^{i, t- i}\ \ (\subseteq E_2^{i, t- i})$}\\
is not necessarily true and so dose not assert that\\
\centerline{$\Ass_R(E_\infty^{i, t- i})\subseteq \Ass_R(E_2^{i, t- i}).$}\\
In the followings, we state some results about finiteness of the set of associated prime ideals of generalized local cohomology modules which, among other things, establish the statements of \emph{\cite [Corollaries 2.4 through 2.7] {Mafi1}}.
\end{rem}



Coung and Hoang in \cite [Theorem 2.4] {CH1} proved the following corollary when $X$ is a finite module.

\begin{cor} \label {3---8}
Let $M$ be a finite $R$--module, $X$ be an arbitrary $R$--module and $n$ be a non-negative integer such that $\lc^{i}_{\fa}(M, X)= 0$ for all $i< n$. Then we have\\
\centerline{$\Ass_R(\lc^{n}_{\fa}(M, X))= \Ass_R(\Ext^n_{R}(M/{\fa}M, X)).$}
\end{cor}

\begin{proof}
Since, by Corollary \ref {2---14}, $\lc^{i}_{\fa+ \Ann_R M}(X)=0$ for all $i<n$,\\
\centerline{$\Hom_R(R/\fa+ \Ann_R M, \lc^{n}_{\fa}(M ,X))\cong \Ext_R^n(M/(\fa+ \Ann_R M) M, X)$}\\
from Theorem \ref {2---10} (b). Thus we have\\
\centerline{$\begin{array}{llll}
\Ass_R(\lc^{n}_{\fa}(M, X))\!\!
&= \ \ \Ass_R(\lc^{n}_{\fa+ \Ann_R M}(M, X))\\
&= \ \ V(\fa+ \Ann_R M)\bigcap \Ass_R(\lc^{n}_{\fa+ \Ann_R M}(M,X))\\
&= \ \ \Ass_R(\Hom_R(R/\fa+ \Ann_R M, \lc^{n}_{\fa}(M, X)))\\
&= \ \ \Ass_R(\Ext_R^n(M/(\fa+ \Ann_R M) M, X))\\
&= \ \ \Ass_R(\Ext^n_{R}(M/{\fa}M,X)),
\end{array}$}\\
as desired.
\end{proof}



In \cite [Theorem 4.5] {CH}, part (a) of the following corollary has been proven when $X$ is a finite module and $R$ is a local ring.

\begin{cor} \label {3---9}
Suppose that $M$ is a finite $R$--module, $X$ is an arbitrary $R$--module and $n$ is a non-negative integer. Assume also that\\
\centerline{$P_n= \displaystyle\bigcup_{i< n} \Supp_R(\lc^{i}_{\fa}(M,X))\ \ (= \displaystyle\bigcup_{i< n}\Supp_R(\Ext^{i}_{R}(M/{\fa}M,X))).$}\\
Then the following statements hold true.
\vspace{-2mm}
\begin{itemize}\setlength{\itemsep}{-2mm}
       \item[\emph{(a)}] $\Ass_R(\lc^{n}_{\fa}(M ,X))\cup P_n= \Ass_R(\Ext^{n}_{R}(M/{\fa}M, X))\cup P_n.$
       \item[\emph{(b)}] $\Ass_R(\lc^{n}_{\fa}(M, X))\subseteq \Ass_R(\Ext^{n}_{R}(M/{\fa}M,X))\cup{P_n}$.
       \item[\emph{(c)}] $\Ass_R(\Ext^{n}_{R}(M/{\fa}M,X))\subseteq \Ass_R(\lc^{n}_{\fa}(M, X))\cup{P_n}$.
       \item[\emph{(d)}] If $\lc^{i}_{\fa}(M,X)$ has finite support for all $i< n$, then $\Ass_R(\lc^{n}_{\fa}(M,X))$ is a finite set if and only if $\Ass_R(\Ext^{n}_{R}(M/{\fa}M,X))$ is a finite set.
\end{itemize}
\end{cor}

\begin{proof}
(a) If $\fp\notin{\underset{i< n}\bigcup \Supp_R(\lc^{i}_{\fa}(M, X))}$, then\\
\centerline{$\Ass_{R_{\fp}}(\lc^{n}_{{\fa}R_{\fp}}(M_{\fp},X_{\fp}))= \Ass_{R_{\fp}}(\Ext^{n}_{R_{\fp}}(M_{\fp}/{(\fa{R_{\fp}})}M_{\fp},X_{\fp}))$}\\
from Lemma \ref {2---5} (b) and Corollary \ref {3---8}. Thus, $\fp$ is not in the left side if and only if it is not in the right side.
\end{proof}



Recall that, an $R$--module $X$ is said to be {\it weakly Laskerian} if the set of associated prime ideals of any quotient module of $X$ is finite (\cite [Definition 2.1]{DM1}). The category of weakly Laskerian $R$--modules is a Serre subcategory of the category of $R$--modules (\cite [Lemma 2.3 (i)] {DM1}) and is denoted by $\mathcal{C}_{w.l}(R)$.

\begin{cor} \label {3---10} \emph{(cf. \cite [Lemma 3.1] {Mafi2})}
Let $M$ be a finite $R$--module, $X$ be an arbitrary $R$--module, and $n$ be a non-negative integer such that $\Ext^{n- i}_{R}(M, \lc^{i}_{\fa}(X))$ is weakly Laskerian for all $i\leqslant n$. Then
\vspace{-2mm}
\begin{itemize}\setlength{\itemsep}{-2mm}
       \item[\emph{(a)}] $\lc^{n}_{\fa}(M, X)$ is weakly Laskerian.
       \item[\emph{(b)}] $\Ass_R(\lc^{n}_{\fa}(M, X))$ is finite.
\end{itemize}
\end{cor}

\begin{proof}
(a) Apply Theorem \ref {2---6} with $\mathcal{S}= \mathcal{C}_{w.l}(R)$.

(b) This is clear from part (a).
\end{proof}



In the following corollary, we generalize \cite [Corollary 2.4] {Mafi1} and \cite [Theorem 3.3] {Mafi2}. Note that, for a finite $R$--modules $M$ and a non-negative integer $n$, $X\in \mathcal{S}$ implies that $\Ext^{n}_{R}(M/\fa M, X)\in \mathcal{S}$ and, by \cite [Proposition 3.4] {HV}, we have $\Ext^{n}_{R}(M/\fa M, X)\in \mathcal{S}$ when $\Ext^{i}_{R}(R/\fa, X)\in \mathcal{S}$ for all $i\leq n$.

\begin{cor} \label {3---11} \emph{(cf. \cite [Corollary 2.4] {Mafi1} and \cite [Theorem 3.3] {Mafi2})}
Let $M$ be a finite $R$--module, $X$ be an arbitrary $R$--module, and $n$ be a non-negative integer such that $\Ext^{n}_{R}(M/\fa M, X)$ and $\lc^{i}_{\fa}(X)$, for all $i< n$, are weakly Laskerian. Then
\vspace{-2mm}
\begin{itemize}\setlength{\itemsep}{-2mm}
       \item[\emph{(a)}] $\Hom_R(R/\fa, \lc^{n}_{\fa}(M, X))$ is weakly Laskerian.
       \item[\emph{(b)}] $\Ass_R(\lc^{n}_{\fa}(M, X))$ is finite.
\end{itemize}
\end{cor}

\begin{proof}
(a) Apply Theorem \ref {2---10} (a)  with $\mathcal{S}= \mathcal{C}_{w.l}(R)$.

(b) Since $\Ass_{R}(\Hom_R(R/\fa, \lc^{n}_{\fa}(M, X)))= V(\fa)\bigcap \Ass_{R}(\lc^{n}_{\fa}(M, X))= \Ass_{R}(\lc^{n}_{\fa}(M, X))$, the assertion follows from part (a).
\end{proof}



\begin{cor} \label {3---12} \emph{(cf. \cite [Corollary 2.5] {Mafi1})}
Suppose that $R$ is a local ring with maximal ideal $\fm$ and $\dim R\leq 2$. Assume also that $M$ is a finite $R$--module and $X$ is an arbitrary $R$--module such that $\Gamma_\fa(X)$ is weakly Laskerian. Then $\Ass_R(\lc^{i}_{\fa}(M, X))$ is finite for all $i$.
\end{cor}

\begin{proof}
This follows from Corollary \ref {3---10} and \cite [Corollaries 2.4 and 2.5] {Mar}.
\end{proof}



\begin{cor} \label {3---13} \emph{(cf. \cite [Corollary 2.6] {Mafi1})}
Suppose that $R$ is a local ring with maximal ideal $\fm$ and $\dim R= n$. Assume also that $M$ is a finite $R$--module and $X$ is an arbitrary $R$--module such that $\lc^j_\fa(X)= 0$ for all $j\neq n- 1, n$. Then $\Ass_R(\lc^{i}_{\fa}(M, X))$ is finite for all $i$.
\end{cor}

\begin{proof}
Follows from Corollary \ref {3---10} and \cite [Corollaries 2.4 and 2.5] {Mar}.
\end{proof}



\begin{cor} \label {3---14} \emph{(cf. \cite [Corollary 2.7] {Mafi1})}
Suppose that $R$ is a local ring with maximal ideal $\fm$ and $\dim_R R/\fa= 1$. Assume also that $M$ is a finite $R$--module and $X$ is an arbitrary $R$--module. Then $\Ass_R(\lc^{i}_{\fa}(M, X))$ is finite for all $i$.
\end{cor}

\begin{proof}
It follows from Corollary \ref {3---10} or Corollary \ref {3---11}.
\end{proof}


\bibliographystyle{amsplain}

\end{document}